\documentclass[letterpaper, 10pt, conference]{ieeeconf}
\IEEEoverridecommandlockouts
\usepackage{graphics} 
\usepackage{epsfig} 
\usepackage{amsmath} 
\usepackage{amssymb}  
\let\labelindent\relax
\usepackage{enumitem}
\usepackage{xcolor}

\newlist{Aenumerate}{enumerate}{1}
\setlist[Aenumerate]{label=C.\arabic*}

\usepackage{amsthm}

\newcommand{\mcl}[1]{\mathcal{ #1}}
\newcommand{\mbf}[1]{\mathbf{ #1}}
\newcommand{\norm}[1]{\ensuremath{\left\Vert #1\right\Vert}}

\newcommand{\bmat}[1]{\begin{bmatrix} #1\end{bmatrix}}
\newcommand{\bmata}[1]{\left\{\begin{matrix}#1\end{matrix}\right\}}

\newcommand{\E}{\mathbb{E}}
\newcommand{\R}{\mathbb{R}}

\newcommand{\Z}{\mathbb{Z}}

\renewcommand{\S}{\mathbb{S}}

\newtheorem{thm}{Theorem}
\newtheorem{defn}{Definition}[thm]
\newtheorem{lem}{Lemma}[thm]

\newtheorem{example}{Example}

\newcommand{\mylen}{3.1pt}

\newcommand{\ip}[2]{\left\langle{#1},{#2}\right\rangle}
\newcommand{\vmatwo}[2]{\begin{bmatrix}
		#1 \\
		#2
\end{bmatrix}}
\newcommand{\vmathree}[3]{\begin{bmatrix}
		#1 \\
		#2 \\
		#3
\end{bmatrix}}
\renewcommand{\th}{\ensuremath{\theta}}
\newcommand{\et}{\eta}
\newcommand{\hinf}{\ensuremath{H_{\infty}}}
\newcommand{\mc}[1]{\ensuremath{\mathcal{#1}}}
\newcommand{\opPzero}{\ensuremath{\mc{P}}}

\newcommand{\opPspq}{\ensuremath{\mc{P}{\tiny\bmata{P,Q,Q^{\top}\\S,R_1,R_2}}}}
\newcommand{\opPpq}[4]{\ensuremath{\mc{P}{\tiny\bmata{#1,#2_1,#2_2\\#4,#3_1,#3_2}}}}

\newcommand{\opPgenpq}[6]{\ensuremath{\mc{P}{\tiny\bmata{#1,#2,#3\\#6,#4,#5}}}}

\newcommand{\barbf}[1]{\bar{\mathbf{#1}}}
\newcommand{\opa}{\bar{\mc{A}}}
\newcommand{\opb}{\bar{\mc{B}}}
\newcommand{\opc}{\bar{\mc{C}}}
\newcommand{\opd}{\bar{\mc{D}}}
\newcommand{\xb}{\barbf{x}}
\newcommand{\yb}{y}
\newcommand{\ub}{w}
\newcommand{\zb}{\ensuremath{\mathbf{z}}}

\newcommand{\myint}{\int_{a}^{b}}
\newcommand{\myinta}[1]{\int_{a}^{#1}}
\newcommand{\myintb}[1]{\int_{#1}^{b}}

\newcommand{\ltwo}{\ensuremath{L_2}}

\allowdisplaybreaks

\title{\LARGE \bf
	A Generalized LMI Formulation for Input-Output Analysis of Linear Systems of ODEs Coupled with PDEs
}

\author{%
	Sachin Shivakumar$^{1}$\hspace{0.5cm}\and
	Amritam Das$^{2}$\hspace{0.5cm}\and
	Siep Weiland$^{2}$\hspace{0.5cm}\and
	\and
	Matthew M. Peet$^{1}$
	\thanks{$^{1}$ Sachin Shivakumar\{sshivak8@asu.edu\} and Matthew M. Peet\{mpeet@asu.edu\} are with School for Engineering of Matter, Transport and Energy, Arizona State University, Tempe, AZ, 85298 USA }
	\thanks{$^{2}$ Amritam Das\{am.das@tue.nl\} and Siep Weiland\{s.weiland@tue.nl\} are with Department of Electrical Engineering, Eindhoven University of Technology, 5600 MB, The Netherlands }
}

\begin{document}

	\maketitle

\begin{abstract}
In this paper, we consider input-output properties of linear systems consisting of PDEs on a finite domain coupled with ODEs through the boundary conditions of the PDE. This framework can be used to represent e.g. a lumped mass fixed to a beam or a system with delay. This work generalizes the sufficiency proof of the KYP Lemma for ODEs to coupled ODE-PDE systems using a recently developed concept of fundamental state and the associated boundary-condition-free representation. The conditions of the generalized KYP are tested using a positive matrix parameterization of bounded operators resulting in a finite-dimensional LMI, feasibility of which implies prima facie provable passivity or $L_2$-gain of the system. No discretization or approximation is involved at any step and there is no conservatism in the theorems. Comparison with other computational methods show that bounds obtained are not conservative in any significant sense and that computational complexity is lower than existing methods involving finite-dimensional projection of PDEs.
\end{abstract}

\section{Introduction}\label{sec:Intro}

Partial Differential Equations (PDEs) model systems with states which vary not only with time, but with respect to some additional independent parameter or parameters. Examples include: beam models, where the states are deflection and rotation; chemical reaction networks, where states are species concentrations; fluid flow, where the state can be velocity or pressure; and time-delay systems, where the state is the history of a finite-dimensional process. By contrast, Ordinary Differential Equations (ODEs) model systems whose states only depend on time, with common examples derived from rigid-body motion or RLC circuits.

Occasionally, we find systems where the dynamics of a PDE with distributed state are coupled to a system of ODEs - often at the boundary of the domain. This can occur naturally, such as in the case of, e.g.; an aircraft, where a rigid fuselage (ODE) is fixed to the root of a flexible wing (PDE); or in a time-delay system, where ODE state feeds directly into the distributed state (history); fluid flow over an accelerating mass~\cite{fluidFlowrotatingcyl}; and heat-exchange devices~\cite{heatexchanger}. In other cases, the ODE-PDE coupling is the result of attempts to design ODE controllers to stabilize a PDE process - such as in the well-developed Port-Hamiltonian framework~\cite{exponentialPHS}. Furthermore, such systems have vector-valued distributed states, representing, e.g. temperature and flow velocities, displacement of flexible structures attached to rigid bodies or concentration of chemical subspecies.

In this paper, we consider a general class of vector-valued linear PDE systems whose internal dynamics are of the form
\begin{align}\label{eq:PDEsimple}
\bmat{\dot{\zb}_1\\\dot{\zb}_2\\\dot{\zb}_3}(s,t) &= A_0(s)\bmat{\zb_1\\\zb_2\\\zb_3}(s,t)+A_1(s)\partial_s\bmat{\zb_{2}\\\zb_{3}}(s,t) \notag\\
&\qquad+ A_2(s) \partial_s^2\zb_{3}(s,t), \qquad \mbf z(s,0)=0;
\end{align}
coupled at the boundary using
{\footnotesize
\begin{align*}
&Bz_b(t) = B_1x(t) \qquad \\
&z_b(t) = \text{col}\big(\mathbf{z}_2(a,t), \mathbf{z}_2(b,t),\mathbf{z}_3(a,t), \mathbf{z}_3(b,t), \partial_s\mathbf{z}_{3}(a,t),\partial_s\mathbf{z}_{3}(b,t))
\end{align*}}
with a linear ODE system
\begin{align*}
	\dot{x}(t) &= Ax(t)+B_2z_b(t), \qquad x(0)=0,
\end{align*}
where $x(t)\in\R^{n_x}$ is the ODE state and $\zb_i(s,t)\in\R^{n_i}$ are the distributed states. Here $A$ is a matrix and $A_0$, $A_1$, $A_2$ are matrix-valued functions. We assume $B$ has row rank $n_2+2n_3$.
The difficulty in analyzing systems of this form is that the boundary conditions (and hence ODE state) do not appear explicitly in the dynamics of the PDE part of the system.

\textbf{Illustrative Example} To illustrate this framework, we use the example of a string coupled with an ODE; see \cite{barreau}. The dynamics of the system can be written as
	\begin{align}\label{eq:ex1}
		&\ddot{w}(s,t) = cw_{ss}(s,t),\notag\\
        &\dot{x}(t) = Ax(t) + Bw(1,t),\notag\\
		& w(0,t) = Kx(t), w_s(1,t) = -c_0\dot{w}(1,t),
	\end{align}
	where $w(s,t)$ is transverse displacement of the string and $x$ is the ODE state.
These equations may be rewritten in the proposed form as
\vspace{-0.25cm}
\begin{align*}
	&\bmat{\dot{\zb}_1\\\dot{\zb}_3}(s,t) = \overbrace{\bmat{0&0\\1&0}}^{A_0(s)}\bmat{\zb_1\\\zb_3}(s,t)+\overbrace{\bmat{c\\0}}^{A_2(s)}\frac{\partial^2}{\partial s^2}\zb_3(s,t)\\
	&\dot{x}(t) = Ax(t)+B\zb_3(1,t),\\
	&\zb_3(0,t) = Kx(t), \zb_{3s}(1,t) = -c_0\zb_1(1,t)
\end{align*} where $\zb_1 = \dot{w}$ and $\zb_3 = w$.

\textbf{Analysis of coupled ODE-PDE Systems}
We propose an algorithm to determine input-output properties of systems of the Form~\eqref{eq:PDEsimple}. Specifically, if we define inputs $w(t)$ and outputs $y(t)$ (See Equation~\ref{eq:originalsys}) then we propose algorithms:
 \begin{enumerate}
 	\item \textbf{$L_2$ Gain:} To find the smallest $\gamma$ such that $\left\Vert y\right\Vert_{\ltwo}\le \gamma \left\Vert w\right\Vert_{\ltwo}$ for all $w \in \ltwo$ and
 	\item \textbf{Passivity:} To verify $\ip{y}{w}_{\ltwo}\ge 0$ for all $w \in \ltwo$, assuming dimensions of inputs and outputs are equal.
 \end{enumerate}



In this paper we present a computationally scalable LMI-based algorithms to prove passivity and obtain a bound on \ltwo ~gain for systems in the form of~\eqref{eq:PDEsimple}. The class of systems includes linear ODEs coupled with the PDEs of parabolic, elliptic and hyperbolic type with Dirichlet, Neumann and mixed boundary conditions. The primary technical difficulty in analysis of coupled ODE-PDE systems is that the coupling occurs through the boundary conditions. Boundary conditions constitute \textit{implicit} constraints on the state of the PDE and are not directly represented in the differential form of Eqn.~\eqref{eq:PDEsimple}. Furthermore, the effect of the PDE state on the ODE state is not bounded as the input to the ODE state is a single unmeasurable point of a larger distributed state. These limitations preclude obvious approaches such as construction of a Lyapunov functionals which depends on the joint ODE-PDE state (See earlier work in~\cite{ahmadi2014,ahmadi2016,gahlawat2017,gahlawat2015}). Our solution to this problem is based on an alternative boundary-condition-free representation of the dynamics using a \text{fundamental state}, $\mbf x_f$, as proposed in~\cite{peet_cdps} (See Section~\ref{sec:fundamentalstate}). In this representation, the effect of the boundary conditions is directly represented in the dynamics of the PDE state, which is defined by a bounded operator. Furthermore, the effect of the boundary on the ODE state is directly represented using a bounded operator on the fundamental state of the PDE. As a result, the KYP and Positive Real (PR) lemmas may be generalized to the dynamics of the ODE-PDE system by abstracting the original matrices to bounded operators(See Section~\ref{sec:IOprop}). 
 
Given our boundary-free representation of the conditions of the KYP and PR lemmas, we parameterize the operator variables in these conditions using the PQRS framework \eqref{def_operator}. This framework uses matrices and matrix-valued polynomials to parameterize operators on $\R^m\times\ltwo^n$ and enforces positivity of such operators using LMI constraints on the matrices and polynomial coefficients. This PQRS framework is used to parameterize both the variables and the system matrices. Then, giving an algorithm for composition and transpose of PQRS operators, we are able to succinctly summarize the resulting LMI conditions in Theorems~\ref{th:sos1} and~\ref{th:sos2}. Numerical examples are then used in Section~\ref{sec:num} to show that the resulting $L_2$-gain bounds are not conservative in any significant sense.

Existing methods for analysis of coupled ODE-PDE systems involve approximating the PDE by a finite-set of ODEs using reduced basis methods, discretization, or other finite-dimensional projection methods \cite{electrocardiogram, trafficflow}. However, approximating a distributed state by a finite set of states often leads to large number of state variables. Furthermore, finite-dimensional projection of PDE may change system properties such as passivity, reachability, observability or stability. Consequently, the properties proven for the ODE approximation of a PDE are not provable for the PDE directly - requiring \textit{a posteriori} tests to bound the errors that occur due to truncation. Methods for bounding the truncation error typically depend on method of discretization and may involve significant conservatism.

Other techniques, such as backstepping method \cite{backstepping1, backstepping2, backstepping3, backstepping4} and sliding mode control method \cite{smc1,smc2}, rely on the use of Lyapunov functionals which are chosen a priori. This limits their utility in optimal controller or observer synthesis.
Other methods, such as \cite{ahmadi2014,ahmadi2016,gahlawat2017,gahlawat2015, baudouin2019stability} use LMIs to search over a given set of Lyapunov functionals. However, these methods are largely restricted to scalar or vector-valued PDE systems with a specific choice of boundary conditions and are not applicable to coupled ODE-PDE systems with a general form of boundary conditions; the LMIs need to be recalculated manually for a different set of boundary conditions.
Finally, \cite{hinfPDESachin} considered input-output analysis of vector-valued linear PDEs (with no ODE coupling). Note that coupled ODE-PDE systems can NOT be posed (even inefficiently) as a special case of the vector-valued PDE framework.

\section{Notation}\label{sec:Notation}
$\mathbb{S}^m\subset \R^{m\times m}$ is the set symmetric matrices. 
For a normed space $X$, define $\ltwo^n[X]$ as the Hilbert space of square integrable $\R^n$-valued functions on $X$ with inner product $\langle x,y \rangle_{\ltwo} = \myint x(s)^{\top} y(s) ds$. The Sobolov spaces are denoted $W^{q,n}[X]:=\{x\in \ltwo^n[X] \mid \frac{\partial^k x}{\partial s^k}\in\ltwo^n[X] \text{ for all }k\le q \}$ with the standard Sobolev inner products.
For a given inner product space, $Z$, the operator$\mc{P}:Z\to Z$ is positive semidefinite (denoted $\mc{P}\succcurlyeq 0$) if $\ip{z}{\mc{P}z}_Z\ge 0$ for all $z\in Z$. Furthermore, we say $\mc{P}:Z\to Z$ is coercive if there exists some $\epsilon>0$ such that $\ip{z}{\opPzero z}_Z \geq \epsilon \Vert z \Vert_Z^2$ for all $z \in Z$. $\mcl L(X,Y)$ is the set of bounded linear operators from $X \rightarrow Y$ where $\mcl L(X)=\mcl L(X,X)$. We use font as a pedagogical aid, where typically $x$ indicates $x \in \R^n$, $\mbf x$ indicates $\mbf x \in L_2^n[X]$ and $\bar{\mbf x}$ indicates $\bar{\mbf x} \in \R^m \times L_2^n[X]$ where $m$ and $n$ are clear from context. Likewise, for operators, $A$ typically indicates $A \in \mcl{L}(\R^n)$ is a matrix, $\mcl A$ indicates $\mcl A \in \mcl{L}(L_2[X])$ and $\bar{\mcl A}$ indicates $\bar{\mcl A} \in \mcl{L}(\R^m \times L_2^n[X])$ (or some variation). The partial derivative $\frac{\partial}{\partial s}\mbf x$ is denoted as $\mbf x_s$.

\section{I/O Properties of an Abstract DPS}\label{sec:IOprop}
Before we present our formulation of the coupled ODE-PDE system, we recall sufficient conditions for Passivity and $L_2$-gain of an abstract Distributed Parameter System (DPS) of the form
\begin{align}\label{eq:abstract}
\bar{\mbf x}_t(t) &= \bar{\mcl A} \bar{\mbf x}(t) + \bar{\mcl B} w(t)\notag\\
y(t) &= \bar{\mcl C} \bar{\mbf x}(t)+D w(t),\qquad \bar{\mbf x}(0)=0,
\end{align}
where, $\bar{\mbf x}(t)\in X$ is the state, $y(t)\in \R^{n_y}$ is the output, $w(t)\in\R^{n_w}$ is the exogenous input to the system, and $\opa:X\to Z, ~\opb:\R^{n_w}\to Z, ~\opc:X\to\R^{n_y}$ and $D\in\R^{n_z \times n_w}$ are linear operators.

\begin{thm} \label{th:KYP}
	Suppose there exists a coercive, self-adjoint linear operator $\mcl P:Z \to Z$ and $\gamma> 0$ such that
\begin{align}\label{eq:loihinf}
  &\ip{z}{\mcl P\bar{\mc{A}}z}_Z+\ip{\bar{\mc{A}}z}{\mcl P z}_Z+\ip{z}{\mcl P\bar{\mc{B}}w}_Z\nonumber\\
  &~+\ip{\bar{\mc{B}}w}{\mcl P z}_Z\le \gamma^2\norm{w}^{2}-\norm{\bar{\mc{C}}z}^2-(\bar{\mc{C}}z)^{\top}(Dw)\nonumber\\
  &~~-(Dw)^{\top}(\bar{\mc{C}}z)-(Dw)^{\top}(Dw)
\end{align} for all $z\in X$ and $w\in\R^{m}$. Then for any $w \in \ltwo^m([0,\infty))$ and $y \in \ltwo^q([0,\infty))$ which satisfy \eqref{eq:abstract} for some $\bar{\mbf x}$, $\Vert y\Vert_{\ltwo} \le \gamma\Vert w\Vert_{\ltwo}$.
\end{thm}
\begin{proof}
	The proof can be found in \cite{hinfPDESachin}.
\end{proof}
\begin{thm} \label{th:KYP2}
 	 Suppose there exists a coercive, self-adjoint linear operator $\mcl P:Z \to Z$ such that
\begin{align}\label{eq:loipass}
&\ip{z}{\mcl P\bar{\mc{A}}z}_Z+\ip{\bar{\mc{A}}z}{\mcl P z}_Z+\ip{z}{\mcl P\bar{\mc{B}}u}_Z+\ip{\bar{\mc{B}}u}{\mcl P z}_Z\nonumber\\
  &~\le(\bar{\mc{C}}z)^{\top}u+u^{\top}(\bar{\mc{C}}z)+(Du)^{\top}u+u^{\top}(Du).
\end{align} for all $z\in X$ and $w\in\R^{m}$. Then for any $w \in \ltwo^m([0,\infty))$ and $y \in \ltwo^m([0,\infty))$ which satisfy \eqref{eq:abstract} for some $\bar{\mbf x}$, $\ip{w}{y}_{L_2}\ge 0$.
\end{thm}
\begin{proof}
	The proof can be found in \cite{hinfPDESachin}.
\end{proof}

\section{Operator Representation of Coupled ODE-PDE Systems}\label{sec:Prelim}

Our generalized formulation for coupled linear ODE-PDE systems consists
of a set of ordinary differential equations and a set of partial differential equations coupled either at the boundary or in-domain. 

For some suitably differentiable functions $x: \R^+ \rightarrow \R^{n_x}$ and $\zb_i: \R\times\R^+ \rightarrow \R^{n_i}$, we represent the internal dynamics of the system, for input $w(t)\in\R^{n_w}$ and output $y(t)\in \R^{n_y}$, as
\begin{align}\label{eq:originalsys}
	&\dot{x}(t) = Ax(t)+ \left(\mc{E}\bmat{\zb_1\\\zb_2\\\zb_3}\right)(t)+B_{wo}w(t)\notag\\
	&\bmat{\dot{\zb}_1\\\dot{\zb}_2\\\dot{\zb}_3}(s,t) = A_0(s)\bmat{\zb_1\\\zb_2\\\zb_3}(s,t)+A_1(s)\bmat{\zb_{2s}\\\zb_{3s}}(s,t) \notag\\
	&\quad+ A_2(s) \zb_{3ss}(s,t)+E(s)x(t)+B_{wp}(s)w(t)\notag\\
	&y(t) = Cx(t)+\left(\mc{C}\bmat{\zb_1\\\zb_2\\\zb_3}\right)(t)+D_ww(t),\notag\\
	&\mc{C}\zb_p :=  C_1 \mathbf{z}_b + \int\limits_{a}^{b} \Bigg(C_a(s) \vmathree{\mathbf{z}_1(s)}{\mathbf{z}_2(s)}{\mathbf{z}_3(s)} + C_b(s) \vmatwo{\mathbf{z}_{2s}(s)}{\mathbf{z}_{3s}(s)}\Bigg) \text{d}s,\notag\\
	&\mc{E}\zb_p :=  E_1 \mathbf{z}_b + \int\limits_{a}^{b} \Bigg(E_a(s) \vmathree{\mathbf{z}_1(s)}{\mathbf{z}_2(s)}{\mathbf{z}_3(s)} + E_b(s) \vmatwo{\mathbf{z}_{2s}(s)}{\mathbf{z}_{3s}(s)}\Bigg) \text{d}s,
\end{align}
where $w\in L_{2}^{n_w}$ and $y\in L_{2}^{n_y}$. For $Z := \ltwo^{n_1}\times W^{1,n_2}\times W^{2,n_3}$ and  $n_z := n_1+n_2+n_3$, the system can be represented by constant matrices $A, C, B_{wo},D_w$ and matrix valued polynomials $A_0,A_1,A_2,B_{wp}, E$ . For linear operators $\mc{E}:Z\to \R^{n_x}$, $\mc{C}:Z\to \R^{n_y}$, $C_a$, $E_a$, $C_b$, $E_b$ are matrix-valued polynomials and $C_1$, $E_1$ are constant matrices of appropriate sizes. This parametric form of $\mc{C}$ of combines boundary valued/distributed output of the state $\zb_p$. 

For the distributed state, $\zb_i$, the boundary conditions are represented in the form
\begin{align*}
Bz_{b}(t) = B_1x(t)+B_{2}w(t)
\end{align*}
where $B,~B_1$ and $B_2$ are matrices of suitable dimensions, $z_{b}=\text{col}\big(\mathbf{z}_2(a), \mathbf{z}_2(b),\mathbf{z}_3(a), \mathbf{z}_3(b), \mathbf{z}_{3s}(a), \mathbf{z}_{3s}(b)\big)$. Note, $B$ should have a row rank of $n_2+2n_3$ for \eqref{eq:originalsys} to have a unique solution. The matrices $B_1$ and $B_2$ represent the coupling with the ODE and disturbance at the boundary, respectively.

We can write this system in the form of \eqref{eq:abstract} as
\begin{align}\label{eq:combined_sys}
\bmat{\dot{x}\\\dot{\zb_p}}(t) &= \opa \bmat{x\\\zb_p}(t) + \opb w(t),\notag\\
y(t) &= \opc \bmat{x\\\zb_p}(t)+\opd w(t),
\end{align}
where $\zb_p = col(\zb_1,\zb_2,\zb_3)$.
Linear operators $\opa:\R^{n_x}\times Z\to \R^{n_x}\times\ltwo^{n_z}$, $\opb:\R^{n_w}\to \R^{n_x}\times\ltwo^{n_z}$, $\opc:\R^{n_x}\times Z\to \R^{n_y}$ and $\opd:\R^{n_w}\to \R^{n_y}$ are defined as
\begin{align}\label{eq:opdef}
	&\left(\opa \bmat{x\\\zb_p}\right)(s) = \bmat{Ax+(\mc{E}\zb_p)(s)\\E(s)x+(\mc{A}\zb_p)(s)},\quad\opd w = D_w w,\notag\\
	&~(\opb w)(s) = \bmat{B_{wo}w\\B_{ws}(s)w} , ~~\quad\left(\opc \bmat{x\\\zb_p}\right) = Cx+(\mc{C}\zb_p),\qquad
\end{align}
where
\begin{align*}
&(\mc{A}\zb_p)(s) := \notag\\
&\hspace{1.2cm} A_0(s)\vmathree{\mathbf{z}_1(s)}{\mathbf{z}_2(s)}{\mathbf{z}_3(s)} + A_1(s)\vmatwo{\mathbf{z}_{2s}(s)}{\mathbf{z}_{3s}(s)}+ A_2(s) \mathbf{z}_{3ss}(s).
\end{align*}



\noindent\textbf{Notation:}
The domain of system \eqref{eq:combined_sys} is
\begin{align}
\label{PDE_boundary}
D_{\opa}:=&\{col(x,\mathbf{z}_p)\in X:  B \ \mathbf{z}_{b} = B_1x(t)+B_{2}w(t) \}
\end{align}
where $X := \mathbb{R}^{n_{x}} \times \ltwo^{n_1}[a, b]\times W^{1,n_2}[a, b]\times W^{2,n_3}[a, b]$. For $\bmat{u\\\mbf{v}},\bmat{x\\\mbf{y}}\in X$, we define the inner product as
\begin{align*}
\Bigg\langle\begin{bmatrix}
x\\
\mbf y
\end{bmatrix}, \begin{bmatrix}
u\\
\mbf v
\end{bmatrix}\Bigg\rangle_{X}:= x^{\top} u + \big\langle \mbf y, \mbf v\big\rangle_{L_2^{n_{z}}[a, b]},
\end{align*}
and norm as
\begin{align*}
\Bigg\lvert\Bigg\lvert\begin{bmatrix}
x\\
\mbf y
\end{bmatrix}\Bigg\rvert\Bigg\rvert_{X}^{2} = \mid\mid x\mid\mid_{\mathbb{R}^{n_{x}}}^2 + \mid\mid \mbf y\mid\mid_{L_2^{n_{z}}[a, b]}^2.
\end{align*}


\begin{example}
	Consider
	\begin{align*}
		&\dot{x}(t) = -x(t)+x(t-\tau)+d(t), ~y(t) = x(t),\\
		&x(s) = 0 \quad ~s\in[-\tau,0].
	\end{align*}
	This can be written as an ODE coupled with a PDE as follows.
	\begin{align*}
	&\dot{x}(t) = -x(t)+z(0,t)+d(t),\quad y(t) = x(t),\\
	&\dot{z}(s,t) = \frac{1}{\tau} z_s(s,t), \quad z(1,t) = x(t), \quad z(s,0) = 0.
	\end{align*}
	This is written in the form of \eqref{eq:PDEsimple} using following matrices.
	\begin{align*}
	&B = [0 ~~1], ~B_1 = 1, ~A = -1, ~B_{wo} = 1,\\
	&A_0(s) = 0, ~A_1(s) =1, ~A_2(s)=0, ~C = 1,E_1 = [1 ~0].
	\end{align*}
	All other matrices are zero.
\end{example}
\section{The PQRS Parametrization of operator}\label{sec:Operator}
Using Theorems \ref{th:KYP} and \ref{th:KYP2}, it is possible to express a bound on \ltwo-gain (or test for passivity) as a test for existence of a coercive, linear operator \opPzero ~that satisfies the corresponding operator inequalities in the Theorem statements. 
We convert operator valued inequality constraints to a set of polynomial equality constraints, which can be converted to LMIs using SOSTOOLS \cite{sostools}, by using matrix valued polynomials to parameterize an operator on $\R\times\ltwo$. This PQRS parametric form is described below. 
\begin{defn}
\label{def_operator}
For a matrix $P \in \mathbb{R}^{n\times m}$, and bounded polynomial functions $Q_1: [a, b] \rightarrow \mathbb{R}^{n\times q}$, $Q_2: [a, b] \rightarrow \mathbb{R}^{r\times m}$, $S: [a, b] \rightarrow \mathbb{R}^{r\times q}$, and $R_1, R_2: [a, b]\times [a, b] \rightarrow \mathbb{R}^{r\times q}$, we define the operator $\opPpq{P}{Q}{R}{S}:\R^{m}\times\ltwo^{q}[a,b]\to \R^{n}\times\ltwo^{r}[a,b]$ as
\begin{align}
    \label{mathcalP}
    &\left(\opPpq{P}{Q}{R}{S}\bmat{x\\\zb}\right)(s):=\\
    & {\small\begin{bmatrix}
    Px + \int_{a}^{b} Q_1(s)\zb(s)ds\\
    Q_2(s)x +S(s)\zb(s)+ \int_{a}^{s} R_1(s,\th)\zb(\th)d\th + \int_{s}^{b} R_2(s,\th)\zb(\th)d\th
    \end{bmatrix}}.\notag
\end{align}
\end{defn}
	
The generalization of KYP and PR Lemmas using abstract DPS form leads to operations, such as composition and adjoint of bounded operators. Hence, we revisit the results related to composition and transpose of PQRS operator; refer \cite{das2018}. These results are used to convert the operator-valued inequalities, in Theorems \ref{th:KYP} and \ref{th:KYP2}, to a PQRS operator positivity constraint that can be enforced using LMIs. First we present sufficient conditions for the positivity of a self-adjoint PQRS operator. Next, we provide the equations to find adjoint of a PQRS operator in Lemma \ref{lem:adjoint} and composition of two PQRS operators in Lemma \ref{lem:composition}. 
\subsection{\protect\opPspq ~Positivity}
\begin{thm}\label{th:positivity}
    For any functions $Z_1:[a, b]\to\R^{d_1\times n}$, $Z_2:[a, b]\times [a, b]\rightarrow\mathbb{R}^{d_2\times n}$, if $g(s)\geq0$ for all $s\in[a,b]$ and 
    \setlength{\abovedisplayskip}{\mylen}
    \setlength{\belowdisplayskip}{\mylen}
   {\small
   \begin{align}
    &P= T_{11}\myint g(s) ds,\nonumber\\
    &Q(\eta) = g(s)T_{12}Z_1(\eta)+\myintb{\eta} g(s)T_{13}Z_2(s,\eta)\text{d}s \notag\\
    &\hspace{1.2cm}+ \myinta{\eta}g(s) T_{14}Z_2(s,\eta)\text{d}s, \nonumber\\
    &R_1(s,\eta) =g(s)Z_1(s)^{\top}T_{23}Z_2(s,\eta)+g(\et)Z_2(\eta,s)^{\top}T_{42}Z_1(\eta)\notag\\
    &\hspace{1.3cm}~+\myintb{s}g(\th)Z_2(\theta,s)^{\top}T_{33}Z_2(\theta,\eta)\text{d}\theta\notag\\
    &\hspace{1.3cm}~~+\int_{\eta}^{s}g(\th)Z_2(\theta,s)^{\top}T_{43}Z_2(\theta,\eta)\text{d}\theta\nonumber\\
    &\hspace{1.3cm}~~~+\myinta{\eta}g(\th)Z_2(\theta,s)^{\top}T_{44}Z_2(\theta,\eta)\text{d}\theta,\nonumber\\
    &R_2(s,\eta) =g(s)Z_1(s)^{\top}T_{32}Z_2(s,\eta)+g(\et)Z_2(\eta,s)^{\top}T_{24}Z_1(\eta)\notag\\
    &\hspace{1.3cm}~+\myintb{\eta}g(\th)Z_2(\theta,s)^{\top}T_{33}Z_2(\theta,\eta)\text{d}\theta\notag\\
    &\hspace{1.3cm}~~+\int_{s}^{\eta}g(\th)Z_2(\theta,s)^{\top}T_{34}Z_2(\theta,\eta)\text{d}\theta\nonumber\\
    &\hspace{1.3cm}~~~+\myinta{s}g(\th)Z_2(\theta,s)^{\top}T_{44}Z_2(\theta,\eta)\text{d}\theta,\nonumber\\
    &S(s) = g(s)Z_1(s)^{\top} T_{22} Z_1(s).
    \label{eq:TH}
   \end{align}}
   where
   \begin{align*}
   T= \begin{bmatrix}
   T_{11} & T_{12} & T_{13} & T_{14}\\
   T_{21} & T_{22} & T_{23} & T_{24}\\
   T_{31} & T_{32} & T_{33} & T_{34}\\
   T_{41} & T_{42} & T_{43} & T_{44}
   \end{bmatrix}\succcurlyeq 0,
   \end{align*}\\
   then the operator \opPspq ~as defined in \eqref{mathcalP} is positive, i.e. $\opPspq\succcurlyeq 0$.
   \end{thm}
   \begin{proof}
   	The proof can be found in \cite{hinfPDESachin}.
   \end{proof}
For convenience, we define the following set.
\begin{align*}
\Phi_d &:= \{{\tiny\bmata{P,Q,Q^{\top}\\S,R_1,R_2}} : \\
&\qquad{\tiny \bmata{P,Q,Q^{\top}\\S,R_1,R_2} = \bmata{P_a,Q_a,Q_a^{\top}\\S_a,R_{1a},R_{2a}}+ \bmata{P_b,Q_b,Q_b^{\top}\\S_b,R_{1b},R_{2b}}}, ~\text{where}\\ &\qquad(P_a,Q_a,S_a,R_{1a},R_{2a}) ~\text{and} ~(P_b,Q_b,S_b,R_{1b},R_{2b}) \\
&\qquad\text{satisfy the conditions of Thm.\ref{th:positivity}}~\text{ with} ~Z_1=Z_d ~\text{and}\\
	 &\qquad\text{where} ~g(s)=1 ~\text{and } g(s) = (s-a)(b-s), \text{resp.}\}
\end{align*}
\subsection{\protect\opPpq{P}{Q}{R}{S} ~composition}
\begin{lem}\label{lem:composition}
Suppose $A,P\in\R^{m\times m}$ are matrices and $B_1, Q_1: [a, b]\to\R^{m\times n}$, $B_2, Q_2: [a, b]\to\R^{n\times m}$, $D,S:[a, b]\to\R^{n \times n}$, $C_i,R_i:[a, b]\times [a, b] \to \R^{n\times n}$ are bounded functions for $i \in \{1,2\}$. Then for any $x \in \mathbb{R}^{m}$ and $\zb \in L_2^n([a, b])$, we have
\begin{align*}
\opPpq{\hat{P}}{\hat{Q}}{\hat{R}}{\hat{S}}\vmatwo{x}{\zb}=
    \opPpq{A}{B}{C}{D}\opPpq{P}{Q}{R}{S}\vmatwo{x}{\zb} 
\end{align*} 
where
\begin{align}\label{eq:composition}
    \hat{P} &= AP + \myint B_1(s)Q_2(s)\text{d}s,\notag\\
    \hat{Q}_1(s) &= AQ_1(s) + B_1(s)S(s)+\myintb{s}B_1(\eta)R_1(\eta,s)\text{d}\eta\notag\\
    &\qquad+\myinta{s}B_1(\eta)R_2(\eta,s)\text{d}\eta,\notag\\
    \hat{Q}_2(s) &= B_2(s)P + D(s)Q_2(s) + \myinta{s}C_1(s,\eta)Q_2(\eta)\text{d}\eta\notag\\
    &\qquad+\myintb{s}C_2(s,\eta)Q_2(\eta)\text{d}\eta,\notag\\
    \hat{S}(s) &= D(s)S(s),\notag\\
    \hat{R}_1(s,\eta) &=B_2(s)Q_1(\eta)+D(s)R_1(s,\eta)+C_1(s,\eta)S(\eta)\notag\\
    &\hspace{-0.5cm}+\myinta{\eta} C_1(s,\theta)R_2(\theta,\eta)\text{d}\theta+\int_{\eta}^{s}C_1(s,\theta)R_1(\theta,\eta)\text{d}\theta\notag\\
    &\hspace{-0.5cm}+\myintb{s}C_2(s,\theta)R_1(\theta,\th)\text{d}\theta,\notag\\
    \hat{R}_2(s,\eta) &=B_2(s)Q_1(\eta)+D(s)R_2(s,\eta)+C_2(s,\eta)S(\eta)\notag\\
    &\hspace{-0.5cm}+\myinta{s} C_1(s,\theta)R_2(\theta,\eta)\text{d}\theta+\int_{s}^{\eta}C_2(s,\theta)R_2(\theta,\eta)d\theta\notag\\
    &\hspace{-0.5cm}+\myintb{\eta}C_2(s,\theta)R_1(\theta,\eta)\text{d}\theta.
\end{align}
\end{lem}
\subsection{\protect\opPpq{P}{Q}{R}{S} ~adjoint}
\begin{lem}\label{lem:adjoint}
Suppose $P\in\R^{m\times m}$ is a matrix and $Q_1: [a, b]\to\R^{m\times n}$, $Q_2: [a, b]\to\R^{n\times m}$, $S:[a, b]\to\R^{n \times n}$, $R_1, R_2:[a, b]\times [a, b] \to \R^{n\times n}$ are bounded functions. Then for any $\barbf{x},\barbf{y} \in \mathbb{R}^{m}\times L_2^n[a, b]$, we have
\begin{align*}
    \Big\langle{\opPpq{\hat{P}}{\hat{Q}}{\hat{R}}{\hat{S}}\barbf{x}},{\barbf{y}}\Big\rangle=\Big\langle{\barbf{x},\opPpq{P}{Q}{R}{S}\barbf{y}}\Big\rangle,
\end{align*}
 where
\begin{align}\label{eq:adjoint}
    &\hat{P} = P^{\top}, &\hat{S}(s) = S(s)^{\top}, \notag\\
    &\hat{Q}_1(s) = Q_2(s)^{\top}, &\hat{R}_1(s,\eta) = R_2(\eta,s)^{\top}, \notag\\
    &\hat{Q}_2(s) = Q_1(s)^{\top}, &\hat{R}_2(s,\eta) = R_1(\eta,s)^{\top}.
\end{align}
\end{lem}
Proof for Lemmas \ref{lem:adjoint} and \ref{lem:composition} can be found in \cite{das2018} \textit{arXiV}  version.

\noindent\textbf{Notation:}
We say,
\begin{enumerate}
	\item $\scriptstyle\bmata{\hat{P},\hat{Q}_1,\hat{Q}_2\\\hat{S},\hat{R}_1,\hat{R}_2}=	\bmata{A,B_1,B_2\\D,C_1,C_2}\times \bmata{P,Q_1,Q_2\\S,R_1,R_2}$ 
	
	if functions $\{\hat{P},\hat{Q}_1,\hat{Q}_2,\hat{S},\hat{R}_1,\hat{R}_2\}$ satisfy \eqref{eq:composition}.
	\vspace{2mm}
	\item $\scriptstyle\bmata{\hat{P},\hat{Q}_1,\hat{Q}_2\\\hat{S},\hat{R}_1,\hat{R}_2}=\bmata{P,Q_1,Q_2\\S,R_1,R_2}^*$ 
	
	if functions $\{\hat{P},\hat{Q}_1,\hat{Q}_2,\hat{S},\hat{R}_1,\hat{R}_2\}$ satisfy \eqref{eq:adjoint}.
\end{enumerate}

\section{Representation of the Dynamics in the Fundamental State}~\label{sec:fundamentalstate}
In this section, we express primal states $\mathbf{z}_p$, as a linear transformation of fundamental states, $\mathbf{z}_f$. The solution to the dynamics of fundamental states always satisfies the boundary conditions of the original PDE. The transformation from fundamental to primal state can be represented as a bounded PQRS operator on $\R^{n_w+n_x}\times\ltwo^{n_z}$. Likewise, the operators $\opa$, $\opb$, $\opc$ and $\opd$ as defined in \eqref{eq:combined_sys} can be represented by some bounded PQRS operator on $\R^{n_w+n_x}\times\ltwo^{n_z}$.

\begin{lem}\label{lem:ftc}
	Suppose $w\in\R^{n_w}$,~$x\in\R^{n_x}$,
	$\text{col}(\zb_1,\zb_{2},\zb_{3})\in \ltwo^{n_1}\times W^{1,n_2}\times W^{2,n_3}$ satisfying
	\begin{align*}
	&B\bmat{\mathbf{z}_2(a)~ \mathbf{z}_2(b)~\mathbf{z}_3(a)~ \mathbf{z}_3(b)~ \mathbf{z}_{3s}(a)~ \mathbf{z}_{3s}(b)}^{\top}\\&\hspace{5cm}=B_1x+B_{2}w
	\end{align*}  where $B$ has a row rank $n_2+2n_3$. Then, for $\mathbf{z}_f=\text{col}(\zb_1,\zb_{2s},\zb_{3ss})$ and $w_r = \text{col}(w,x),$
	\begin{align*}
	&\bmat{\zb_1\\\zb_{2}\\\zb_{3}} = \opPgenpq{0}{0}{H_0}{G_1}{G_2}{G_0}\bmat{w_r\\\zb_f},\\
	&\qquad\bmat{\zb_{2s}\\\zb_{3s}} = \opPgenpq{0}{0}{H_1}{G_4}{G_5}{G_3}\bmat{w_r\\\zb_f},
	\end{align*}	
	where
	\begin{align}\label{eq:FTC}
	&H_0(s) = K(s)(BT)^{-1}\bmat{B_{2}&B_1},\notag\\
	&H_1(s) = V(s)(BT)^{-1}\bmat{B_{2}&B_1},\notag\\
	&G_0(s) = \bmat{I & 0 & 0\\0 &0&0\\0 &0&0},~ G_3(s) = \bmat{0 & I & 0\\0 &0&0},\notag\\
	&G_2(s,\th) = -K(s)(BT)^{-1}BQ(s,\th),\notag\\
	&G_5(s,\th) = -V(s)(BT)^{-1}BQ(s,\th),\notag\\
	&G_1(s,\th) = \bmat{0 & 0 & 0\\0 &I&0\\0 &0&(s-\th)I}+G_2(s,\th),\notag\\
	&G_4(s,\th) = \bmat{0 & 0 & 0\\0 &0&I}+G_5(s,\th),\notag\\
	&K(s) = \bmat{0 & 0 & 0\\I &0&0\\0 &I&(s-a)I},\quad V(s) = \bmat{0 & 0 & 0\\0 &0&I}\notag\\
	&T = \bmat{I &0&0\\I &0&0\\0 & I & 0\\0 &I&(b-a)I\\0&0&I\\0&0&I},\quad Q(s,\th) = \bmat{0 &0&0\\0 &I&0\\0 & 0 & 0\\0 &0&(b-\th)I\\0&0&0\\0&0&I}.
	\end{align}
\end{lem}
Lemma \ref{lem:ftc} can be proved by using fundamental theorem of calculus; refer \cite{peet_cdps}. Now, clearly, for any $\zb_f\in L_2^{n_1+n_2+n_3}[a,b]$, $\zb_p\in L_2^{n_1}[a,b]\times W^{1,n_2}\times W^{2,n_3}$ and it satisfies the boundary conditions. Rewriting the dynamics in terms of $\zb_f$, essentially, eliminates the need for boundary conditions.5

\opPspq ~notation can also be used to rewrite the system operators, as defined in \eqref{eq:opdef}, in a compact form.
\begin{lem}\label{lem:opnotation}
	Suppose the operators $\opa, \opb, \opc, \opd$ are as defined in \eqref{eq:opdef}. For any $w\in\R^{n_w},~x\in\R^{n_x},\text{col}(\zb_1,\zb_{2},\zb_{3})\in\ltwo^{n_1}\times W^{1,n_2}\times W^{2,n_3} $ such that
	\begin{align*}
	&B\bmat{\mathbf{z}_2(a)~ \mathbf{z}_2(b)~\mathbf{z}_3(a)~ \mathbf{z}_3(b)~ \mathbf{z}_{3s}(a)~ \mathbf{z}_{3s}(b)}^{\top}\\&\hspace{5cm}=B_1x+B_{2}w
	\end{align*}  where $B$ has a row rank $n_2+2n_3$, if
	\begin{align*}
	&\zb_f=\bmat{\zb_1\\\zb_{2s}\\\zb_{3ss}},~\zb_p=\bmat{\zb_1\\\zb_{2}\\\zb_{3}}, ~w_r = \bmat{w\\x},
	\end{align*}
	then
	\begin{align}\label{eq:opNotation}
	&\opa \bmat{x\\\zb_p}  = \opPgenpq{\hat{A}_0}{\hat{A}_1}{\hat{A}_2}{\hat{A}_4}{\hat{A}_5}{\hat{A}_3}\bmat{w_r\\\zb_f},~\opb \ub= \opPgenpq{\hat{B}_0}{0}{\hat{B}_2}{0}{0}{0}\bmat{w_r\\\zb_f},\nonumber\\
	&\opc \bmat{x\\\zb_p}=   \opPgenpq{\hat{C}_0}{\hat{C}_1}{0}{0}{0}{0}\bmat{w_r\\\zb_f},~\opd \ub= \opPgenpq{\hat{D}_0}{0}{0}{0}{0}{0}\bmat{w_r\\\zb_f},\notag\\
	&w = \opPgenpq{I_0}{0}{0}{0}{0}{0} \vmatwo{w_r}{\zb_{f}},~\bmat{x\\\zb_p}=\opPgenpq{I_1}{0}{H_0}{G_1}{G_2}{G_0}\bmat{w_r\\\zb_f},
	\end{align}
	where
	\begin{align}\label{eq:sysop}
	&\scriptstyle\bmata{\hat{A}_0,\hat{A}_1,\hat{A}_2\\\hat{A}_3,\hat{A}_4,\hat{A}_5}\notag\\&= \scriptstyle\bmata{0,E_a,0\\A_0,0,0}\times\bmata{0,0,H_0\\G_0,G_1,G_2}+\bmata{0,E_b,0\\A_1,0,0}\times\bmata{0,0,H_1\\G_3,G_4,G_5}\notag\\
	&\qquad\scriptstyle+\bmata{E_1,0,0\\0,0,0}\times\bmata{0,0,0\\T_1,T_2,T_2}+\bmata{[0 ~A],0,[0 ~E]\\ [0 ~0 ~A_2],0,0},\notag\\
	&\scriptstyle\bmata{\hat{C}_0,\hat{C}_1,0\\0,0,0}\notag\\
	&= \scriptstyle\bmata{0,C_a,0\\0,0,0}\times\bmata{0,0,H_0\\G_0,G_1,G_2}+\bmata{0,C_b,0\\0,0,0}\times\bmata{0,0,H_1\\G_3,G_4,G_5}\notag\\
	&\qquad\scriptstyle+\bmata{C_1,0,0\\0,0,0}\times\bmata{0,0,0\\T_1,T_2,T_2}+\bmata{[0 ~C],0,0\\ 0,0,0},\notag\\
	&\hat{B}_0 = \bmat{B_{wo}&0},~\hat{B}_2(s) = \bmat{B_{wp}(s)&0},~\hat{D}_0 = \bmat{D_w & 0},\notag\\
	&I_0= \bmat{I &0}, ~ I_1= \bmat{0&I},~T_1 = T\bmat{B_2 &B_1},\notag\\
	&T_2(s,\th) = T(BT)^{-1}BQ(s,\th)+Q(s,\th).\notag\\
	\end{align}
\end{lem}
\begin{proof}
	This can be proved using Lemma \ref{lem:ftc}
\end{proof}

\noindent\textbf{Notation:}
For convenience, we define the following short-hand notation, for the operators in Lemma \ref{lem:opnotation}.
\begin{align}\label{eq:shorthand}
	&\mc{P}^A:=\opPgenpq{\hat{A}_0}{\hat{A}_1}{\hat{A}_2}{\hat{A}_4}{\hat{A}_5}{\hat{A}_3},\quad&\mc{P}^B:=\opPgenpq{\hat{B}_0}{0}{\hat{B}_2}{0}{0}{0},\nonumber\\
	&\mc{P}^C:=\opPgenpq{\hat{C}_0}{\hat{C}_1}{0}{0}{0}{0},&\mc{P}^D:=\opPgenpq{\hat{D}_0}{0}{0}{0}{0}{0},\notag\\
	&\mc{P}^I:=\opPgenpq{I_0}{0}{0}{0}{0}{0},&\mc{P}^0:=\opPgenpq{I_1}{0}{H_0}{G_1}{G_2}{G_0},
\end{align}

\noindent\textbf{Example.} For the example in section \ref{sec:Prelim},
\begin{align*}
&\dot{x}(t) = -\frac{1}{\tau}\int_{0}^{1} z_s(s,t) ds+d(t),\quad y(t) = x(t),\\
&\dot{z}(s,t) = \frac{1}{\tau} z_s(s,t), \quad z(1,t) = x(t), \quad z(s,0) = 0,
\end{align*}
we can find the matrices defined in Lemma \ref{lem:opnotation}.
\begin{align*}
&\hat{A}_0 = \bmat{0&0},~\hat{A}_1(s) = -\frac{1}{\tau},~\hat{A}_2(s) = 0,\hat{A}_3(s) = \frac{1}{\tau},\notag\\
&~\hat{A}_4(s,\th) = 0,~\hat{A}_5(s,\th) = 0,\hat{B}_0 = \bmat{1&0},\notag\\
&\hat{B}_2(s) = \bmat{0&0},~\hat{C}_0 = \bmat{0&1},~\hat{C}_1(s)=0,~\hat{D}_0 = \bmat{0 & 0}.
\end{align*}
For this system, the boundary constraint $z(1,t) = x(t)$ is directly embedded in the dynamics via $\hat{A}_0$ and $\hat{A}_1$.

\section{Reformulation of the LOI for a coupled ODE-PDE}\label{sec:reform}
Using the PQRS parametric form, the sufficient conditions in Theorems \ref{th:KYP} and \ref{th:KYP2}, which are in the form of operator feasibility tests, can be posed as polynomial constraints. Then, these constraints can be used to test for passivity or find a bound on \ltwo ~gain of the system \eqref{eq:combined_sys}.

\begin{thm}\label{th:sos1}
	Suppose there exists $\epsilon>0$, $\gamma>0$, $d_1, d_2 \in \Z$, matrix $P\in\S^{n_x}$, matrix-valued polynomials $Q:\R\to\R^{n_x\times n_z}$, $S:\R\to\S^{n_z}$, and $R_1,R_2:\R\times\R\to\R^{n_z\times n_z}$ such that
	\begin{align*}
	{\bmata{P-\epsilon I,Q,Q^{\top}\\S-\epsilon I,R_1,R_2}} \in \Phi_{d_1}.\end{align*}
	Then for all $\xb(t)\in X$, $y\in \ltwo^{n_y}[0,\infty)$ and $\ub\in \ltwo^{n_w}[0,\infty)$ which satisfy \eqref{eq:combined_sys}, if
	\begin{align*}
		-\bmata{J_0,J_1,J_2\\J_3,J_4,J_5}-\bmata{J_0,J_1,J_2\\J_3,J_4,J_5}^*&\in \Phi_{d_2},
	\end{align*}
	where
	\begin{align*}
	&\scriptstyle\bmata{J_0,J_1,J_2\\J_3,J_4,J_5}= \scriptstyle\bmata{K_0,K_1,K_2\\K_3,K_4,K_5}\times\bmata{\hat{A}_0+\hat{B}_0,\hat{A}_1,\hat{A}_2+\hat{B}_2\\\hat{A}_3,\hat{A}_4,\hat{A}_5}\\
	&\qquad+\frac{1}{2}\scriptstyle\bigg(\bmata{\hat{C}_0^{\top}\hat{C}_0,0,\hat{C}_1^{\top}\hat{C}_0\\0,0,0}+\bmata{\hat{D}_0^{\top}\hat{D}_0-\gamma^2I_0,0,0\\0,0,0}\bigg),\\
	&~\qquad+\scriptstyle\bmata{\hat{C}_0,\hat{C}_1,0\\0,0,0}^*\times\bmata{\hat{D}_0,0,0\\0,0,0},\\
	&\bmata{K_0,K_1,K_2\\K_3,K_4,K_5}= \bmata{I_1,0,H_0\\G_0,G_1,G_2}^*\times\bmata{P,Q,Q^{\top}\\S,R_1,R_2}
	\end{align*}
	$H_i, G_i$ are as defined in Lemma \ref{lem:ftc}, $\hat{A}_i, \hat{B}_i, \hat{C}_i, \hat{D}_i$ and $I_i$ are as defined in Lemma \ref{lem:opnotation}, then $\Vert y\Vert_{\ltwo} \le \gamma\Vert \ub\Vert_{\ltwo}$.
\end{thm}
\begin{proof} The proof is in the Appendix.
\end{proof}
\begin{thm}\label{th:sos2}
	Suppose there exists $\epsilon>0$, $\gamma>0$, $d_1, d_2 \in \Z$ matrix $P\in\S^{n_x}$, matrix-valued polynomials $Q:\R\to\R^{n_x\times n_z}$, $S:\R\to\S^{n_z}$, and $R_1,R_2:\R\times\R\to\R^{n_z\times n_z}$ such that
	\begin{align*}
	{\bmata{P-\epsilon I,Q,Q^{\top}\\S-\epsilon I,R_1,R_2}} \in \Phi_{d_1}.\end{align*}
	Then for all $\xb(t)\in X$, $y\in \ltwo^{n_y}[0,\infty)$ and $\ub\in \ltwo^{n_y}[0,\infty)$ which satisfy \eqref{eq:combined_sys}, if
	\begin{align*}
	-\bmata{J_0,J_1,J_2\\J_3,J_4,J_5}-\bmata{J_0,J_1,J_2\\J_3,J_4,J_5}^*&\in \Phi_{d_2},
	\end{align*}
	where
	\begin{align*}
	&\scriptstyle\bmata{J_0,J_1,J_2\\J_3,J_4,J_5}=\scriptstyle\bmata{K_0,K_1,K_2\\K_3,K_4,K_5}\times\bmata{\hat{A}_0+\hat{B}_0,\hat{A}_1,\hat{A}_2+\hat{B}_2\\\hat{A}_3,\hat{A}_4,\hat{A}_5}\\
	&\hspace{2.5cm}-\scriptstyle\bmata{(\hat{C}_0+\hat{D}_0)^{\top}I_0,0,\hat{C}_1^{\top}I_0\\0,0,0},\\
	&\scriptstyle\bmata{K_0,K_1,K_2\\K_3,K_4,K_5}= \scriptstyle\bmata{I_1,0,H_0\\G_0,G_1,G_2}^*\times\bmata{P,Q,Q^{\top}\\S,R_1,R_2},
	\end{align*}
	$H_i, G_i$ are as defined in \eqref{eq:FTC}, $\hat{A}_i, \hat{B}_i, \hat{C}_i, \hat{D}_i$ and $I_i$ are as defined in \eqref{eq:sysop}, then $\ip{\ub}{y}_{\ltwo}\geq0$.
\end{thm}
\begin{proof}
	The proof is similar to that of Theorem \ref{th:sos1}, see Appendix.
\end{proof}
\section{Numerical Implementation and Testing}\label{sec:num}
In this section, we validate the accuracy of proposed algorithm. It was implemented in MATLAB using an adaptation of SOSTOOLS \cite{sostools}. The code can be found on CodeOcean (\textit{https://codeocean.com/capsule/4730069/}).
\subsection{Stabilizing Boundary Control of PDEs}
Consider the example from \cite{boundarycontrol} stabilized by a backstepping controller at the boundary. The resulting closed-loop system is as follows.
\begin{align*}
&\dot{x}(t) = -3x+w(0,t)+d(t),\\
&\dot{w}(s,t) = w_{ss}(s,t)+d(t),\\
&w_s(0,t) = 0, ~w(L,t) = 0,\\
&y(t) = \int_{0}^{L} w(s,t) ds
\end{align*}
A bound on \ltwo ~gain of this closed-loop system in presence of disturbance $d(t)$, using the proposed method, was found to be 0.4269 for relatively low order monomial basis ($d_1$=2). On the other hand, the norm bound obtained through a finite difference method (approximately 100 discrete elements) had significant conservatism with a value of 0.5941.
\subsection{Aircraft wings as flexible beams}

Consider a simplified model of an aircraft, in which a lumped mass is attached to flexible Euler-Bernoulli beams on either side. The deflection of beam can be used to estimate the stresses that develop in the beam. Let $z$ represent vertical displacement of the aircraft, $w$ be the deflection of the wings on either side and $d(t),~u(t)$ disturbances.
\begin{align*}
&\ddot{z}(t) =-EIw_{sss}(0,t)+ d(t), \\
&\ddot{w}(s,t) = - \frac{EI}{\mu} w_{ssss} + u(t),\\
& w(0,t) = z(t), w_s(0,t) = 0, w_{ss}(L,t)=0, w_{sss}(L,t) =0\\
&y(t) = w(L,t)
\end{align*}
where we assume displacement of the wings are symmetric.
The use of variables, $x_1=z$, $x_2=\dot{z}$, $v_1=\dot{w}$ and $v_2=w_{ss}$, converts the system into
\begin{align*}\small\setlength\arraycolsep{0.25pt}
&\vmatwo{\dot{x}_1}{\dot{x}_2}(t) = \bmat{0&1\\0&0}\vmatwo{x_1}{x_2}(t)-\bmat{0&0\\0&F}\int_{0}^{L}v_{2ss}(s,t)ds\\
&\hspace{3cm}+\bmat{0\\d(t)},\\
&\vmatwo{\dot{v}_1}{\dot{v}_2}(s,t) = \bmat{0&-\frac{EI}{\mu}\\1&0}\vmatwo{v_{1ss}}{v_{2ss}}(s,t)+\bmat{u(t)\\0}\\
&y(t)=\int_{0}^{L} (L-s)v_2(s,t)ds.
\end{align*}
This is a linear system of ODEs coupled with PDEs at the boundary for which a bound on \ltwo ~gain can be found using proposed framework. Using tip displacement as the output, we can estimate the stresses tip in presence of disturbances. For $\frac{EI}{\mu} = 10$, the \ltwo ~gain was found to be 0.8936.
\subsection{Time-delay systems}
We test conservatism of the bounds by comparing our \hinf ~norm bound to the method described in \cite{peettactds} for a few well studied time-delay systems. We use the following three examples to document the result from numerical tests in Table I.

\begin{enumerate}[label=C.\arabic*:]
	\item 	
	$\dot{x}(t) = -x(t)-x(t-\tau)+d(t), ~y(t) = x(t).$
	
	\item 
	{\small\begin{flalign*}
	&\dot{x}(t) = \bmat{0&1\\-2&0.1}x(t) + \bmat{0&0\\1&0}x(t-\tau)+\bmat{1&0\\0&1}w(t),\\
	&y(t) = \bmat{0&1}x(t).
	\end{flalign*}}
	\item
	{\small\setlength\arraycolsep{1.25pt}
	\begin{flalign*}
	&\dot{x}(t) = \bmat{-2&0\\0&-0.9}x(t) + \bmat{-1&0\\-1&-1}x(t-\tau)+\bmat{-0.5\\1}w(t),\\
	&y(t) = \bmat{1&0}x(t).
	\end{flalign*}}
\end{enumerate}
\vspace{-0.25cm}
\begin{table}[!h]\label{tab:tds}
	\begin{center}
	\begin{tabular}{|c|c|c|}
		\hline
		          &ODE-PDE framework   & Method in \cite{peettactds} \\\hline
		C.1&0.8911&0.8920\\\hline
		C.2&2.9366&2.9367\\\hline
		C.3&0.2601&0.2601\\\hline
	\end{tabular}
	\caption{Comparison of bounds obtained by different methods for Time-delay systems at delay value, $\tau=1$.}
	\end{center}
\end{table}
\section{Conclusion}\label{sec:conclu}
In this paper, we propose a method to prove passivity and obtain bounds for the \ltwo -gain of coupled linear ODE-PDE systems with disturbances at the boundary or in-domain distributed using the LMI framework. The method presented does not use discretization and the properties established are prima facie provable. Restricting the PQRS operator to polynomial basis can, in theory, lead to some amount conservatism. However, the numerical results indicate the bounds are not conservative in any significant sense.
\section*{Acknowledgments}
This work was supported by Office of Naval Research Award N00014-17-1-2117 and National Science Foundation under grant No. 1739990.
\bibliographystyle{ieeetr}
\bibliography{references}
\appendix
The procedure of proof for Theorem \ref{th:sos1}, stated in Section \ref{sec:reform}, is described below.
\begin{proof}
	From Theorem \ref{th:KYP}, if we can find a self-adjoint, coercive operator $\opPzero$ such that
	\begin{align*}
	&\ip{\xb}{\opPzero\opa\xb}_X+\ip{\opa\xb}{\opPzero \xb}_X+\ip{\xb}{\opPzero\opb\ub}_X+\ip{\opb\ub}{\opPzero\xb}_X\nonumber\\
	&~~-( \gamma^2\ub^{\top}\ub-(\opc\xb)^{\top}(\opc\xb)-(\opc\xb)^{\top}(\opd\ub)-(\opd\ub)^{\top}(\opc\xb)\nonumber\\
	&~~~~-(\opd\ub)^{\top}(\opd\ub) )\leq 0
	\end{align*} for all $\xb\in X$ and $\ub\in \R^{n_w}$, then $\norm{y}\leq \gamma\norm{w}$.
	Let $\opPzero = \opPspq$.
	Then
	\begin{align*}
	&\ip{\bmat{x\\\zb_p}}{\opPspq\bmat{x\\\zb_p}}\\
	&= \ip{\bmat{x\\\zb_p}}{\opPgenpq{P-\epsilon I}{Q}{Q^{\top}}{R_1}{R_2}{S-\epsilon I}\bmat{x\\\zb_p}}+\epsilon \norm{x}^2+\epsilon\norm{\zb_p}^2\\
	&\ge \epsilon \norm{x}^2+\epsilon\norm{\zb_p}_{\ltwo}^2.
	\end{align*}
	\opPspq ~is self-adjoint and coercive.
	Next, we prove that
	\begin{align*}
	&\ip{\xb}{\opPspq\opa\xb}_X+\ip{\opa\xb}{\opPspq \xb}_X\nonumber\\
	&~+\ip{\xb}{\opPspq\opb\ub}_X+\ip{\opb\ub}{\opPspq\xb}_X\nonumber\\
	&~~-( \gamma^2\ub^{\top}\ub-(\opc\xb)^{\top}(\opc\xb)-(\opc\xb)^{\top}(\opd\ub)-(\opd\ub)^{\top}(\opc\xb)\nonumber\\
	&~~~~-(\opd\ub)^{\top}(\opd\ub) \leq 0
	\end{align*} for all $w(t)\in\R^{n_w}$ and $\xb(t)\in X$.
	
	Let $\barbf{x}_f = \text{col}(w,x,\zb_1,\zb_{2s},\zb_{3ss})$. Then from Lemma \ref{lem:opnotation} and short-hand defined in \eqref{eq:shorthand},
	\begin{align*}
	&\ip{\xb}{\opPspq\opa\xb}_X =\ip{\mc{P}^0\barbf{x}_f}{\opPspq\mc{P}^A\barbf{x}_f}\\
	&\hspace{3.3cm}=\ip{\barbf{x}_f}{(\mc{P}^0)^*\opPspq\mc{P}^A\barbf{x}_f}.\\
	&\text{Similarly,} \\
	&\ip{\xb}{\opPspq\opb\ub}_X =\ip{\mc{P}^0\barbf{x}_f}{\opPspq\mc{P}^B\barbf{x}_f}\\
	&\hspace{3.3cm}=\ip{\barbf{x}_f}{(\mc{P}^0)^*\opPspq\mc{P}^B \barbf{x}_f},\\
	&\norm{\opc\xb}^2= \left(\mc{P}^C\barbf{x}_f\right)^{\top}\left(\mc{P}^C\barbf{x}_f\right)=\ip{\barbf{x}_f}{(\mc{P}^C)^*\mc{P}^C\barbf{x}_f},\\
	&\left(\opc\xb\right)^{\top}\left(\opd\xb\right)= \left(\mc{P}^C\barbf{x}_f\right)^{\top}\left(\mc{P}^D\barbf{x}_f\right)=\ip{\barbf{x}_f}{(\mc{P}^C)^*\mc{P}^D\barbf{x}_f},\\
	&(\opd\xb)^{\top}(\opd\xb)= \left(\mc{P}^D\barbf{x}_f\right)^{\top}\left(\mc{P}^D\barbf{x}_f\right)=\ip{\barbf{x}_f}{(\mc{P}^D)^*\mc{P}^D\barbf{x}_f},\\
	&\text{and}\\
	&(\ub)^{\top}(\ub)= \left(\mc{P}^I\barbf{x}_f\right)^{\top}\left(\mc{P}^I\barbf{x}_f\right)=\ip{\barbf{x}_f}{(\mc{P}^I)^*\mc{P}^I\barbf{x}_f}.
	\end{align*}
	Then
	\begin{align*}
	&\ip{\xb}{\opPspq\opa\xb}_X+\ip{\opa\xb}{\opPspq \xb}_X\nonumber\\
	&~+\ip{\xb}{\opPspq\opb\ub}_X+\ip{\opb\ub}{\opPspq\xb}_X\nonumber\\
	&~~-( \gamma^2\ub^{\top}\ub-(\opc\xb)^{\top}(\opc\xb)-(\opc\xb)^{\top}(\opd\ub)-(\opd\ub)^{\top}(\opc\xb)\nonumber\\
	&~~~~-(\opd\ub)^{\top}(\opd\ub))
	=\ip{\barbf{x}_f}{\left((\mc{P}_{eq})^{*}+\mc{P}_{eq}\right)\barbf{x}_f}_{\ltwo} 		
	\end{align*}
	where $\mc{P}_{eq} = \opPgenpq{J_0}{J_1}{J_2}{J_4}{J_5}{J_3}$.
	From the Theorem statement, $\mc{P}_{eq}+\mc{P}_{eq}^*\preccurlyeq0$.
	
	Then, all the conditions of Theorem \ref{th:KYP} are satisfied. Hence $\Vert \yb\Vert_{\ltwo} \le \gamma\Vert \ub\Vert_{\ltwo}$.
\end{proof}
\end{document}